\documentclass{article} 
\usepackage{amssymb}
\usepackage{amsmath}
\begin{document} 
\input epsf.sty

\title{On polynomials sharing 
preimages of compact sets,
and related questions
}
 
\author{F. Pakovich}
\date{}

\maketitle

\def\be{\begin{equation}}
\def\ee{\end{equation}}
\def\bs{$\square$ \vskip 0.2cm}
\def\d{{\rm d}} 
\def\D{{\rm D}} 
\def\I{{\rm I}} 
\def\C{{\mathbb C}} 
\def\N{{\mathbb N}} 
\def\P{{\mathbb P}}
\def\Z{{\mathbb Z}}
\def\R{{\mathbb R}} 
\def\ord{{\rm ord}}
\def\ssigma{\omega}

\def\e{\eqref}
\def\phi{{\varphi}}
\def\v{{\varepsilon}} 
\def\deg{{\rm deg\,}} 
\def\Det{{\rm Det}}
\def\dim{{\rm dim\,}} 
\def\Ker{{\rm Ker\,}} 
\def\Gal{{\rm Gal\,}}
\def\St{{\rm St\,}} 
\def\exp{{\rm exp\,}} 
\def\cos{{\rm cos\,}} 
\def\diag{{\rm diag\,}} 
\def\GCD{{\rm GCD }}
\def\LCM{{\rm LCM }}
\def\mod{{\rm mod\ }}

\def\bp{\begin{proposition}}
\def\ep{\end{proposition}}
\newtheorem{zzz}{Theorem}
\newtheorem{yyy}{Corollary}
\def\bt{\begin{theorem}}
\def\et{\end{theorem}}
\def\be{\begin{equation}}
\def\bee{\begin{equation*}}
\def\la{\label}
\def\l{\lambda}
\def\m{\mu}
\def\ee{\end{equation}}
\def\eee{\end{equation*}}
\def\bl{\begin{lemma}}
\def\el{\end{lemma}}
\def\bc{\begin{corollary}}
\def\ec{\end{corollary}}
\def\pr{\noindent{\it Proof. }}
\def\note{\noindent{\bf Note. }}
\def\bd{\begin{definition}}
\def\ed{\end{definition}}

\newtheorem{theorem}{Theorem}[section]
\newtheorem{lemma}{Lemma}[section]
\newtheorem{definition}{Definition}[section]
\newtheorem{corollary}{Corollary}[section]
\newtheorem{proposition}{Proposition}[section]

\begin{abstract}
In this paper we give a solution of
the following problem: under what conditions on 
infinite compact sets $K_1,K_2\subset \C$ and polynomials
$f_1,$ $f_2$ the preimages $f_1^{-1}\{K_1\}$ and $f_2^{-1}\{K_2\}$
coincide.
Besides, we investigate some related questions. In particular, we show
that polynomials sharing an invariant compact set 
distinct from a point have equal Julia sets. 
\end{abstract}

\section{Introduction} Let $f_1(z),$ $f_2(z)$ be complex polynomials
and $K_1,K_2\subset \C$ be finite or infinite compact sets. 
In this paper we investigate  
the following problem.
Under what conditions on the collection $f_1(z), f_2(z), K_1, K_2$ 
the preimages $f_1^{-1}\{K_1\}$ and $f_2^{-1}\{K_2\}$
coincide that is \be \la{1} f_1^{-1}\{K_1\}=f_2^{-1}\{K_2\} =K\ee
for some compact set $K\subset \C$ ? 
Let us mention several particular cases 
when the answer is known.

The following problem was posed 
in \cite{y}: whether the equality
$f_1^{-1}\{-1,1\}=f_2^{-1}\{-1,1\}$ for polynomials of the same degree 
$f_1(z),f_2(z)$ implies that
$f_1(z)=\pm f_2(z)$ ? This problem was solved in \cite{p1}, \cite{p2}.
It was shown that actually for any compact set $K\subset \C$ containing at 
least $2$ points and polynomials of the same degree $f_1(z),f_2(z)$
the equality $f_1^{-1}\{K\}=f_2^{-1}\{K\}$ 
implies that $f_1(z)=\sigma(f_2(z))$ for some linear function $\sigma(z)=az+b,$ $a,b\in \C,$ 
such that $\sigma\{K\}=K.$

For polynomials of arbitrary degrees solutions of the equation
$f_1^{-1}\{K\}=f_2^{-1}\{K\}$ for a compact set $K\subset \C$
of the positive logarithmic capacity
were described in \cite{d}. Recently this result was extended 
to an arbitrary infinite compact set $K$ in \cite{d2}.
It was shown that if $K$ is distinct from a union of
circles or a segment and $\deg f_2(z) \geq \deg f_1(z)$ then there
exists a polynomial $g(z)$ such that $f_2(z)=g(f_1(z))$ and 
$g^{-1}\{K\}=K.$

Furthermore, the problem of description of pairs of
polynomials $f_1(z),$ $f_2(z)$ sharing the Julia set, 
studied in \cite{b}, \cite{f}, \cite{be1}, \cite{be2}, \cite{sh},
\cite{a2}, \cite{a1}, also is a particular case of 
problem \eqref{1}. Here the answer (\cite{sh},
\cite{a2}) says that whenever the common Julia set $J$ is distinct from 
a circle or a segment there exists a polynomial 
$p(z)$ such that $J$ is the Julia set of $p(z)$ and
up to a symmetry of $J$ the 
polynomials $f_1(z)$ and $f_2(z)$ are the iterations of $p(z).$  

Finally, notice that problem (1)
absorbs the classical problem of description of commuting polynomials  
(\cite{j}, \cite{fa}, \cite{r2}, \cite{er}) 
since commuting polynomials are known to have equal Julia sets.

In this paper we provide a surprisingly simple description of solutions of equation \eqref{1}
which in particular permits to treat and reprove all the results
mentioned above in 
the uniform way. Namely, we relate equation \eqref{1} to the functional equation
\be \la{2}
g_1(f_1(z))=g_2(f_2(z)), \ee
where $f_1(z),f_2(z),g_1(z),g_2(z)$ are polynomials. 
It is easy to see that for any polynomial
solution of \eqref{2} and any compact set $K_3\subset \C$
we obtain a solution of \eqref{1} setting \be \la{3} K_1=g_1^{-1}\{K_3\}, \ \ \
K_2=g_2^{-1}\{K_3\}. \ee In particular, for any ``decomposable" polynomial $f_2(z)=g_1(f_1(z))$ and any compact set $S\subset \C$
we have: 
$$f_2^{-1}\{S\}=f_1^{-1}\{T\},$$ where $T=g_1^{-1}\{S\}.$ 
  
The main result of this paper
states that, under a very mild condition
on the cardinality of $K,$
all solutions of \eqref{1} can be obtained in this way.
Moreover, 
using the Ritt theory of factorisation of polynomials
we describe
these solutions in a very explicit way.

\begin{zzz} Let $f_1(z),$ $f_2(z)$ be
polynomials, $\deg f_1=d_1,$ $\deg f_2=d_2,$ \linebreak $d_1\leq d_2,$ and $K_1,K_2, K\subset \C$ be compact sets
such that \eqref{1} holds. 
Suppose that
$\rm{card}\{K\}\geq \LCM(d_1,d_2).$ Then, if $d_1$ divides $d_2,$
there exists a polynomial $g_1(z)$ such that $f_2(z)=g_1(f_1(z))$
and $K_1=g_1^{-1}\{K_2\}.$ On the other hand, if $d_1$ does not divide $d_2,$ then there exist polynomials $g_1(z),$ $g_2(z),$
$\deg g_1=d_2/d,$ $\deg g_2=d_1/d,$ where $d=\GCD(d_1,d_2),$ and a compact set $K_3\subset \C$
such that \eqref{2},\eqref{3} hold. Furthermore, in this case there
exist polynomials
$\tilde f_1(z),$ $\tilde f_2(z),$ $W(z),$ $\deg W(z)=d,$
such that 
\be \la{kom} f_1(z)=\tilde f_1(W(z)), \ \ \ \ f_2(z)=\tilde f_2(W(z))\ \ \ \   
\ee 
and 
there exist linear functions $\sigma_1(z),$ $\sigma_2(z)$ such that 
either
\begin{align} \la{t2}  & g_1(z)=z^cR^{d_1/d}(z) \circ \sigma_1^{-1},
&& \tilde f_1(z)  =\sigma_1 \circ z^{d_1/d}, \\
& g_2(z)=z^{d_1/d} \circ \sigma_2^{-1},  && \tilde f_2(z)  =\sigma_2 \circ z^cR(z^{d_1/d}), 
\notag \end{align}
for some 
polynomial $R(z)$ and $c$ 
equal to the remainder after division of $d_2/d$ by $d_1/d,$
or
\begin{align} \la{t3} \ \ \ \ \ \  & g_1(z)= T_{d_2/d}(z)\circ \sigma_1^{-1}, \ \ \ \ \ \ \ \ \ && \tilde f_1(z)  
=\sigma_1 \circ T_{d_1/d}(z),  &&   \\
& g_2(z)= T_{d_1/d}(z)\circ \sigma_2^{-1}, \ \ \ \ \ \ \ \ \ &&  \tilde f_2(z)  =\sigma_2 \circ T_{d_2/d}(z),
\notag\end{align}
for the Chebyshev polynomials $T_{d_1/d}(z),$ $T_{d_2/d}(z).$

\end{zzz}

As a corollary of theorem 1 we obtain the following simple
description of the solutions
of \eqref{1} with $d_1=d_2$ (cf. (\cite{p2}, \cite{p1}). In particular,
this description implies the results of \cite{f}, \cite{be1}
concerning the polynomials of the same degree sharing the Julia set. 

\begin{yyy} If
equality \eqref{1}
holds for
polynomials $f_1(z),$ $f_2(z)$ such that $d_1=d_2$
and at least one of the sets $K_1$, $K_2$ 
contains more than one point then 
there exists a linear function $\sigma(z)$
such that $f_2(z)=\sigma(f_1(z))$ and $K_2=\sigma\{K_1\}.$
\end{yyy}

As an other corollary of theorem 1 we describe the situations when the preimage of a compact set under a polynomial mapping may have symmetries. This result generalizes the corresponding results of \cite{b}, \cite{be2} proved under assumption that $K$ is the Julia set of $f(z).$

Denote by $\Sigma_T$ the group
of linear functions which transform the set $T\subset \C$ to itself. 

\begin{yyy} Let $f(z)$ be a polynomial and $K, K_1\subset \C$ be compact sets such that $K=f^{-1}\{K_1\}.$     
Then $\Sigma_{K}$ is a group of rotations. Furthermore,  either $K$ 
is a union of circles and $f(z)=\sigma_2 \circ z^{d_1} \circ \sigma_1$
for some linear functions $\sigma_1(z)$, $\sigma_2(z),$
or $\Sigma_{K}$ is finite and $f(z)=\sigma_2 \circ z^aR(z^b) \circ \sigma_1$ 
for some linear functions $\sigma_1(z)$, $\sigma_2(z)$ and a polynomial
$R(z),$ where $b$ equals the order of $\Sigma_{K}$ and $a<b.$

\end{yyy}

In the case 
when $K_1=K_2$ in \eqref{1} the totality
of solutions of the corresponding equation
\be \la{u22} f_1^{-1}\{T\}=f_2^{-1}\{T\}=K \ee
becomes much smaller in comparison with the general case.
Namely, under notation
introduced above the following result holds.

\begin{zzz} 
Let $f_1(z),$ $f_2(z)$ be
polynomials such that \eqref{u22} holds for some infinite compact sets 
$T, K\subset \C.$ Then, if $d_1$ divides $d_2,$
there exists a polynomial
$g_1(z)$ such that $f_2(z)=g_1(f_1(z))$ and $g_1^{-1}\{T\}=T.$
On the other hand, if $d_1$ does not divide $d_2,$ then
there exist polynomials
$\tilde f_1(z),$ $\tilde f_2(z),$ $W(z),$ $\deg W(z)=d,$ 
satisfying \eqref{kom}. 
Furthermore, in this case one of the following conditions holds.
\vskip 0.2cm \noindent
1) $T$ is a union of circles with the common center and 
\be \la{tt2} \tilde f_1(z)=\sigma \circ z^{d_1/d}, \ \ \ \ \ \tilde f_2(z)=\sigma \circ \gamma z^{d_2/d}
\ee for some linear function $\sigma(z)$ and $\gamma \in \C.$ 
\vskip 0.2cm \noindent
2) $T$ is a segment and 
\be \la{tt3} \tilde f_1(z)=\sigma \circ \pm T_{d_1/d}(z), \ \ \ \ \ \tilde f_2(z)=\sigma \circ \pm T_{d_2/d}(z),\ee
for some linear function $\sigma(z)$
and the Chebyshev polynomials $T_{d_1/d}(z),T_{d_2/d}(z).$
\vskip 0.2cm
\end{zzz}

This result was also obtained in \cite {d}, \cite{d2}. 
However, our method is completely different from the method 
used in these papers. In particular, in our proof we do not use the
classification of commuting polynomials 
that eventually allows us to obtain   
a new proof of  
this classification.

Furthermore, we 
describe the polynomials sharing an invariant compact set
that is solutions of the equation \be \la{xx} f_1^{-1}\{T\}=f_2^{-1}\{T\}=T. \ee 
where $f_1(z), f_2(z)$ are polynomials
and $T\subset \C$ is any compact set. 

The theorem 3 below generalizes 
results of \cite{sh},
\cite{a2} proved under the assumption that 
$T$ is the Julia set of $f_1(z), f_2(z).$

\begin{zzz}  Let $f_1(z), f_2(z)$ be polynomials and $T\subset \C$ 
be a compact set such that \eqref{xx} holds.  
Then  
one of the following conditions holds.
\vskip 0.2cm
\noindent 1) $T$ is a union of circles and
\be \la{po1} f_1(z)=\sigma \circ z^{d_1}\circ \sigma^{-1}, \ \ \ \ \ f_2(z)=\sigma \circ \gamma z^{d_2} \circ \sigma^{-1} 
\ee for some linear function $\sigma(z)$ and $\gamma \in \C,$ where $\vert
\gamma\vert=1$ whenever $T$ is distinct from a point.
\vskip 0.2cm 
\noindent 2) 
$T$ is a segment and
\be \la{po2} f_1(z)=\sigma \circ \pm T_{d_1}\circ \sigma^{-1}, \ \ \ \ \ f_2(z)=\sigma \circ \pm T_{d_2} \circ \sigma^{-1} 
\ee for some linear function
$\sigma(z)$ and the Chebyshev polynomials $T_{d_1}(z),T_{d_2}(z).$ 
\vskip 0.2cm
\noindent 3) The group $\Sigma_T$ is finite and 
there exist a polynomial $p(z)$ and integers $s_1,s_2$ such that $p^{-1}\{T\}=T$ and
\be \la{po3}
f_1(z)= \mu_1 \circ p^{\circ s_1} , \ \ \ f_2(z)=\mu_2 \circ p^{\circ s_2}
\ee
for some linear functions $\mu_1(z), \mu_2(z)\in\Sigma_T.$   
\end{zzz}

It was shown in \cite{b}, \cite{be2} that 
polynomials sharing the Julia set are closely related to the functional equation \be \la{urr}  f_1(f_2(z))=\mu(f_2(f_1(z))), \ee
where $\mu(z)$ is a linear function. It turns out that the same is true for equation \eqref{xx}. Furthermore, theorem 4 below states that 
actually polynomials sharing an invariant compact set 
have the same Julia sets and that any of these properties is equivalent to equation
\eqref{urr} for an appropriate linear function $\mu(z).$ 
Note that together with theorem 3 this implies in particular 
the classification 
of commuting polynomials (cf. \cite{j}, \cite{fa}, \cite{r2},
\cite{er}). 
\begin{zzz} 
The following conditions are equivalent:
\vskip 0.1cm 
\noindent 1) Equality \eqref{xx} holds for some
compact set $T\subset \C$ distinct from a point. 
\vskip 0.1cm 
\noindent 2) 
Polynomials $f_1(z), f_2(z)$ have the same Julia sets.
\vskip 0.1cm 
\noindent 3) There exist compact sets $T_1,T_2 \subset \C$ such that 
$f_1^{-1}\{T_1\}=T_1,$ $f_2^{-1}\{T_2\}=T_2$ and equation \eqref{urr} holds
for some $\mu(z)\in \Sigma_{T_1}\cap \Sigma_{T_2}$.
\end{zzz}

The approach of this paper is similar to the one introduced by the
author for solving the Yang problem cited above. It consists in using 
a relation between a polynomial $f(z)$ and the $n$-th polynomial of least deviation $p_n(z)$ on the preimage $f^{-1}\{K\}$ of a compact set $K\subset \C$ (see section 2 below). 
This relation 
together with the uniqueness theorem for the $n$-th polynomial of
least deviation and the Ritt theorem
permits to reduce equation \eqref{1} to equation \eqref{2}.

The paper is organized as follows. In the second section 
we recall some classical results about polynomial approximations
and prove theorem 2.3 which generalizes the previous results of papers \cite{peh}, \cite{kamo}, \cite{p2}, \cite{p1}. Although essentially we need only the weaker previous result from \cite{kamo} we give the proof 
of theorem 2.3 
because we believe that this result is interesting by itself. 

In the third section we recall two theorems about polynomial solutions of
equation \eqref{2} which we use subsequently. 
In the fourth section using approach
described above we give the proofs of theorem 1 and corollaries 1,2.

In the fifth section we prove theorem 2. Here, the idea of the proof is to examine
the infinite chain of compact sets and polynomials obtained by repeated use of theorem 1. 

Finally, in the sixth section using the obtained results as well as some constructions from 
the papers on Julia sets cited above we prove theorems 3, 4.

\section{Polynomial approximations}
Denote by $P_n$ the vector space consisting of polynomials of degrees 
$\leq n.$ It is known (see e.g. \cite{c}) that for any compact set
$R\subset \C$ and any complex-valued function $\phi(z)$ continuous on $R$ there exists a polynomial
$p_{n,\phi}(z)\in P_n$ such that \be \la{ap}\parallel \phi-p_{n,\phi} \parallel
=\min\limits_{p\in P_n}\parallel \phi-p \parallel,\ee where 
the symbol $\parallel g \parallel$ denotes the uniform norm of 
the function $g$ on $R:$
$$\parallel g \parallel=\max\limits_{x\in R}\vert g(z)\vert.$$
Such a polynomial is called the $n$-th polynomial of least deviation from $\phi$ 
on $R.$ The $n$-th polynomial of least deviation from $\phi$
is known to be unique whenever 
$R$ contains at least $n+1$ points (see e.g. \cite{c}). 
In case when $\phi(z)=z^n$
the polynomial $z^n-p_{n-1,\phi}(z)$ is called the $n$-th monic polynomial of least 
deviation from zero on $R.$ 

It turns out that for an arbitrary compact set $R$
any polynomial $P(z)$ is the polynomial of least deviation 
on the set $P^{-1}\{R\}$ whenever $R$ is ``centered"
at the origin. 
More precisely, the following theorem
holds.

\bt \la{op} Suppose that $R\subset \C$ is a compact set such that the disk of 
the smallest radius 
which contains $R$ is 
centered at the origin. Then any monic polynomial $P(z)$ of degree
$n$ is the $n$-th monic polynomial of least deviation from zero on the set $P^{-1}\{R\}.$
\et

This theorem was proved in \cite{p2} where it was applied 
to the description of solutions of \eqref{u22} with $d_1=d_2.$ 

Note that theorem 2.1 implies the following well known result:
the $n$-th normalized 
Chebyshev polynomial $T_n(z)$ is the $n$-th monic polynomial 
of least deviation from 
zero on $[-1,1].$ Indeed, it is enough to observe that the
formula $T_n(\cos z)=\cos nz$ implies that
$T_n^{-1}[-1,1]=[-1,1].$ 
Similarly, one can deduce that the polynomial $z^n$ is the $n$-th monic polynomial 
of least deviation from zero on any union of circles centered at the origin. 

A more general than theorem 2.1 result was proved by a different method 
(actually, earlier) 
in \cite{kamo} 
in connection with the description of polynomials of least deviation on Julia sets. 

\bt \la{kb} Let $R\subset \C$ be a compact set and $T(z)$ be 
the $m$-th monic polynomial of least deviation from zero on $R.$ Then for any 
polynomial $P(z)$ of degree $n$ with leading coefficient $c_n$ the polynomial 
$T(P(z))/c_n^m$ is the $mn$-th monic polynomial of least deviation from zero 
on the set $P^{-1}\{R\}.$
\et

Finally, some more general result - theorem \ref{mt} below - 
was proved in \cite{peh}. Nevertheless, the proof was given only under 
the additional 
assumption that the so called extremal signature (see e.g. \cite{c}) for $\phi(z)-p_{m,\phi}(z)$ on $R$ contains 
no critical values of $P(z).$  
Below, we give the proof in the general case generalizing the method of \cite{p2}.

\bt \la{mt} Let $R\subset \C$ be a compact set,
$\phi(z)$ be a continuous function on $R,$ and $p_{m,\phi}(z)$ be 
the $m$-th polynomial of least deviation from $\phi(z)$ on $R.$ Then for any 
polynomial $P(z)$ of degree $n$ the polynomial $p_{m,\phi}(P(z))$ is the $mn+n-1$-th 
polynomial of least deviation from $\phi(P(z))$ on the set $P^{-1}\{R\}.$
\et
\noindent{\it Proof of theorem 2.3.}  For any polynomial $Q(z)$ set  
$$Q_P(z)=\frac{1}{n}\sum_{\substack {y\in \C, \\ P(y)=P(z)}}Q(y),$$
where the root $y$ of multiplicity $k$ of $P(y)-P(z)=0$ is repeated $k$ times.

Clearly,  
$$\max\limits_{z\in P^{-1}\{R\}}\vert  \phi(P(z))
-Q_P(z)\vert\leq 
\max\limits_{z\in P^{-1}\{R\}}\sum_{\substack {y\in \C, \\ P(y)=P(z)}}
\frac{\vert  \phi(P(y))
-Q(y)\vert}{n}
.$$
On the other hand, since $P^{-1}\{R\}$ together with a point $z$ contains all
the points $y$ such that $P(y)=P(z),$ 
we have: $$
\max\limits_{z\in P^{-1}\{R\}}\sum_{\substack {y\in \C, \\ P(y)=P(z)}}
\frac{\vert \phi(P(y))
-Q(y)\vert}{n}
\leq
\max\limits_{z\in P^{-1}\{R\}}\vert  \phi(P(z))-Q(z)\vert.
$$
Therefore, for any $Q(z)$ the inequality
\be \la{1a}
\max\limits_{z\in P^{-1}\{R\}}\vert  \phi(P(z))
-Q_P(z)\vert\leq \max\limits_{z\in P^{-1}\{R\}}\vert  \phi(P(z))-Q(z)\vert
\ee holds.

Furthermore, observe that for any polynomial $R(z)$ of degree $<n$
the function $R_P(z)$ is constant. 
Indeed, for $R(z)=z^j,$ $0\leq j \leq n-1,$ this follows from the Newton
formulas which express $R_P(z)$ via the symmetric functions 
$S_j,$ $0\leq j \leq n-1,$ of roots $y_i,$ $1\leq i \leq n,$ of
$P(y)-P(z)=0$ and in general case by the linearity.

Let now $Q(z)$ be a polynomial of an arbitrary degree $q$
and let $$Q(z)=\sum_{i=0}^{[q/n]}a_i(z)P^i(z)$$ be its $P$-adic
decomposition. Then 
$$Q_P(z)=
\sum_{i=1}^{[q/n]}a_iP^i(z),$$ where $a_i=a_{i_P}(z)/n$
are constants. Therefore, 
\be\la{2a} \max\limits_{z\in P^{-1}\{R\}}\vert  \phi(P(z))-
\sum Q_P(z)\vert=
\max\limits_{z\in R}\vert  \phi(z)-\sum_{i=1}^{[q/n]}a_iz^i\vert
.\ee

Suppose now that $q< mn +n.$ Then 
\be \la{3a} \max\limits_{z\in R}\vert  \phi(z)-\sum_{i=1}^{[q/n]}a_iz^i\vert
\geq \max\limits_{z\in R}\vert  \phi(z)-p_{m,\phi}(z)\vert=
\max\limits_{z\in P^{-1}\{R\}}\vert \phi(P(z))-p_{m,\phi}(P(z))\vert.\ee
It follows now from \eqref{1a},\eqref{2a} and \eqref{3a}
that 
$$
\max\limits_{z\in P^{-1}\{R\}}\vert  \phi(P(z))-Q(z)\vert \geq
\max\limits_{z\in P^{-1}\{R\}}\vert \phi(P(z))-p_{m,\phi}(P(z))\vert.
$$
\vskip 0.4cm

\noindent{\it Proof of theorem 2.2.}  Theorem \ref{kb} follows 
from theorem \ref{mt}. Indeed,
for $h(z)\in P_m$ we have:
$$
\min\limits_{p\in P_m}\parallel (\phi+h)-p \parallel=
\min\limits_{p\in P_m}\parallel (\phi+h)-(p+h) \parallel
= \min\limits_{p\in P_m}\parallel \phi-p \parallel
$$
and $p_{m,\phi+h}(z)=p_{m,\phi}(z)+h(z).$ Similarly, for $\beta \in \C$ we have:
$$
\min\limits_{p\in P_m}\parallel \beta \phi-p \parallel=
\min\limits_{p\in P_m}\parallel \beta \phi- \beta p \parallel
= \beta \min\limits_{p\in P_m}\parallel \phi-p \parallel
$$
and $p_{m,\beta \phi}(z)=\beta p_{m,\phi}(z).$ 

Therefore, 
$$p_{mn-1,z^{mn}}(z)=\frac{p_{mn-1,c_n^mz^{mn}}(z)}{c_n^m}=
\frac{p_{mn-1,P^m+ c_n^mz^{mn}-P^m}(z)}{c_n^m}= $$ $$=
\frac{p_{mn-1,P^m}(z)}{c_n^m}+z^{mn}-\frac{P^m(z)}{c_n^m} 
=\frac{p_{m-1,z^m}(P(z))}{c_n^m}+z^{mn}-\frac{P^m(z)}{c_n^m}.
$$ 
Hence, $$z^{mn}-p_{mn-1,z^{mn}}(z)=
\frac{(z^m-p_{m-1,z^m}(z))}{c_n^m}\circ P(z)
=\frac{T(P(z))}{c_n^m}.$$

\vskip 0.4cm

\noindent{\it Proof of theorem 2.1.} Theorem
\ref{op} is a particular case of theorem \ref{kb}
since its condition is equivalent to the condition that 
the first monic polynomial of least deviation from zero on $R$ is $z.$

\section{Solutions of $A(B(z))=C(D(z))$}

In this section we recall two theorems about polynomial 
solutions of the equation 
\be \la{r} A(B(z))=C(D(z)) \ee
proved in \cite{en}, \cite{ri} (see
also \cite{sch}, Theorems 5 and 8).

\bt \la{eng}  
Let $A(z),B(z),C(z),D(z)$ be polynomials such that \eqref{r} holds. Then there exist polynomials
$V(z), \hat B(z), \hat D(z), $ such that 
$$B(z)=\hat B(V(z)), \ \  D(z)=\hat D(V(z)), \ \  \deg V(z)=\GCD(\deg B(z),\deg D(z)),$$ 
and there exist
polynomials $U(z), \hat A(z), \hat C(z)$ such that 
$$A(z)=U(\hat A(z)), \ \  C(z)=U(\hat C(z)), \ \  \deg U(z)=\GCD(\deg A(z),\deg C(z)).$$
\et

The theorem \ref{eng} reduces the problem of finding the solutions of \eqref{r}
to the one when 
$\deg A(z)=\deg D(z)$ and
$\deg B(z)=\deg C(z)$ are coprime. The answer to the last question
is given by the following ``second Ritt theorem".

\bt 
Let $A(z),B(z),C(z),D(z)$ be non-linear
polynomials satisfying  \eqref{r} such that 
$a=\deg A(z)=\deg D(z)$ and $b=\deg B(z)=\deg C(z)$
are coprime and $a>b.$ Then there exist linear functions $\sigma_1(z),\sigma_2(z),\mu(z), \nu(z)$ 
such that either
$$ A(z)=\nu \circ z^cR^b(z) \circ \sigma_1^{-1}, \ \ \
B(z)=\sigma_1 \circ z^b \circ \mu, $$
\be \ \ \  \ \ \    C(z)=\nu \circ z^b \circ \sigma_2^{-1}, \ \ \ \ \ \ \ \ \ \ D(z)=\sigma_2 \circ
z^cR(z^b) \circ \mu \ee 
for some polynomial $R(z)$ and $c$ equal to the remainder after
division of $a$ by $b$, or
$$ A(z)=\nu \circ T_a(z) \circ \sigma_1^{-1}, \ \ \
B(z)=\sigma_1 \circ T_b(z) \circ \mu, $$
\be C(z)=\nu \circ T_b(z) \circ \sigma_2^{-1}, \ \ \ D(z)=\sigma_2 \circ
T_a(z) \circ \mu \ee
for the Chebyshev polynomials $T_a(z), T_b(z).$
\et

\section{Solutions of $f_1^{-1}\{K_1\}=f_2^{-1}\{K_2\} =K$} 
{\it Proof of theorem 1.} Let $ p_{1}(z)$ be the 
$d_2/d$-th monic polynomial of least deviation from zero on $K_1$ and $ p_{2}(z)$ 
be the $d_1/d$-th monic polynomial of least deviation from zero on $K_2.$
Then by theorem \ref{kb} the polynomial $ p_{1}(f_1(z))/a_1^{d_2/d},$ 
where $a_1$ is the 
leading coefficient of $f_1(z),$ is the $d_1d_2/d$-th monic polynomial
of least deviation from zero on $K.$
Similarly, the polynomial $ p_{2}(f_2(z))/a_2^{d_1/d},$ where $a_2$ is
the leading coefficient of $f_1(z),$ is the $d_1d_2/d$-th monic polynomial
of least deviation from zero on $K.$ Since 
$$\rm{card}\{K\}\geq \,\LCM(d_1,d_2)=d_1d_2/d$$ 
it follows from the uniqueness of the polynomial of least deviation that 
\be \la{l} \hat g_{1}(f_1(z))=\hat g_{2}(f_2(z)),\ee 
where
$\hat g_1(z)=p_{1}(z)/a_1^{d_2/d},$ $\hat g_2(z)= p_{2}(z)/a_2^{d_1/d}.$ 
Hence, by theorem 3.1 there exist polynomials
$\tilde f_1(z), \tilde f_2(z), V(z)$ such that 
$$f_1(z)=\tilde f_1(V(z)), \ \ \ f_2(z)=\tilde f_2(V(z)), 
$$ where $\deg V(z)=d.$ 

If $d_1$ divides $d_2$ then the polynomial $\tilde f_1(z)$ is linear and
setting $g_1(z)=\tilde f_2\circ \tilde f_1^{-1}$ we see that
$f_2(z)=g_1(f_1(z)).$ Moreover, since for any polynomial $f(z)$ and sets $T_1,T_2 \subset \C$ the equality $f^{-1}\{T_1\}=f^{-1}\{T_2\}$ implies that $T_1=T_2$, it follows from
the equality 
$$f^{-1}_2\{K_2\}=f^{-1}_1\{g_1^{-1}\{K_2\}\}=f^{-1}_1\{K_1\}=K$$ that 
$K_1=g_1^{-1}\{K_2\}.$

Furthermore, if $d_1$ does not divide $d_2$ then both $\tilde f_1(z), \tilde f_2(z)$ are non-linear and therefore
$\hat g_1(z), \hat g_2(z)$ are also non-linear. Since equality \eqref{l} 
implies the equality 
\be \la{gg} \hat g_1(\tilde f_1(z))=\hat g_2(\tilde f_2(z)),\ee where
$\deg \hat g_1(z)=\deg \tilde f_2(z)$ and $\deg \hat g_2(z)=\deg \tilde f_1(z)$ 
are coprime,
applying theorem 3.2 to \eqref{gg} and 
setting
$$g_1(z)= \nu^{-1}\circ \hat g_1(z), \ \ \  g_2(z)= \nu^{-1}\circ
\hat g_2(z), \ \ \ 
W(z)=\mu \circ V $$ we see that \eqref{2} and \eqref{kom} hold with
$\tilde f_1(z), \tilde f_2(z), g_1(z), g_2(z)$ satisfying either
\eqref{t2} or \eqref{t3}.

Observe now that 
$$ g_1\{K_1\}= g_1\{ f_1\{K\}\}= g_2\{ f_2\{K\}\}=
g_2\{K_2\}.$$ 
Set $K_3=g_1\{K_1\}=g_2\{K_2\}$ and show 
that the equalities 
$$g_1^{-1}\{K_3\}=K_1, \ \ \ g_2^{-1}\{K_3\}=K_2$$
hold.
Notice that it is enough to prove only one of these 
equalities. Indeed, \eqref{2} implies that
\be \la{zzz}
f_1^{-1}\{g_1^{-1}\{K_3\}\}=f_2^{-1}\{g_2^{-1}\{K_3\}\}.\ee
Therefore, if say $g_1^{-1}\{K_3\}=K_1$ then \eqref{1} and
\eqref{zzz}
imply that $$K=f_2^{-1}\{g_2^{-1}\{K_3\}\}.$$ 
Since 
$K=f_2^{-1}\{K_2\}$ 
it follows that  $g_2^{-1}\{K_3\}=K_2.$

Show first that if \eqref{t2} holds then
\be \la{p} g_2^{-1}\{K_3\}=K_2. \ee
Clearly, equality \eqref{p} is equivalent to the equality
$$
\sigma_2^{-1}\{K_2\}=(z^{d_1/d})^{-1}\{K_3\}.
$$ On the other hand,
the last equality
is equivalent to the statement that the set $\sigma_2^{-1}\{K_2\}$ together with a point $x$
contains any point of the form $\varepsilon x,$ 
where $\varepsilon$ is a $d_1/d$-th root of unity.

In order to prove the last statement first observe that
in view of \eqref{1} and \eqref{kom}
we have:
$$
W^{-1}\{\tilde f_1^{-1}\{K_1\}\}=W^{-1}\{\tilde
f_2^{-1}\{K_2\}\}.$$ Hence
$$
\tilde f_1^{-1}\{K_1\}=\tilde
f_2^{-1}\{K_2\}=W\{K\}$$ or equivalently 
\be \la{gh1} (z^{d_1/d})^{-1}\{\sigma_1^{-1}\{K_1\}\}=
(z^cR(z^{d_1/d}))^{-1}\{\sigma_2^{-1}\{K_2\}\}=W\{K\}.\ee
Suppose now that $x\in \sigma_2^{-1}\{K_2\}$ and let $y$ be a point of 
$W\{K\}$ such that $y^cR(y^{d_1/d})=x.$ Then equality \eqref{gh1} implies
that any point of the form $\varepsilon y,$ where $\varepsilon$ is a $d_1/d$-th root of unity, 
also belongs to $W\{K\}.$ Since
$$\sigma_2^{-1}\{K_2\}
=z^cR(z^{d_1/d})\{W\{K\}\} 
$$ it follows that $\sigma_2^{-1}\{K_2\}$ together with a point $x$
contains any point of the form 
$\varepsilon^c y^cR(y^{d_1/d})=\varepsilon^c x.$
To finish the proof it is enough to observe that the equality $\GCD(d_1/d,d_2/d)=1$ 
implies the equality $\GCD(c,d_1/d)=1.$ Therefore, if $\varepsilon$ runs all $d_1/d$-th roots of unity
then $\varepsilon^c$ also runs all $d_1/d$-th roots of unity. 

In the case when \eqref{t3} holds the proof of the equality
$$\sigma_2^{-1}\{K_2\}=(T_{d_1/d})^{-1}\{K_3\}$$
which in this case is equivalent to equality \eqref{p} is similar. We must show that for any point 
$x\in \sigma_2^{-1}\{K_2\}$ all the points
$y$ such that $T_{d_1/d}(y)=T_{d_1/d}(x)$ also belong to $\sigma_2^{-1}\{K_2\}.$ Equivalently, we must
show that if $\cos\, \alpha =x\in \sigma_2^{-1}\{K_2\}$ for some
$\alpha\in \C$ then for any $k=1,2, ... , (d_1/d)-1$ the number 
$\cos (\alpha+\frac{2\pi d}{d_1}k)$ also belongs to $\sigma_2^{-1}\{K_2\}.$ 

As above observe that 
\be \la{gh11} (T_{d_1/d})^{-1}\{\sigma_1^{-1}\{K_1\}\}=
(T_{d_2/d})^{-1}\{\sigma_2^{-1}\{K_2\}\}=W\{K\}.\ee
Suppose now that 
$\cos\, \alpha =x\in \sigma_2^{-1}\{K_2\}$ and set $t=\cos (
\frac{\alpha d}{d_2}).$ 
Then $T_{d_2/d}(t)=x$ and hence $t\in W\{K\}.$  
Therefore, 
\eqref{gh11} 
implies that all the points of the form $$\cos (\frac{\alpha d}{d_2}+ \frac{2\pi d}{d_1}j), \ \ \ j=1,2, ... ,\frac{d_1}{d}-1,$$
belong to $W\{K\}.$ It follows now from the equality
$$\sigma_2^{-1}\{K_2\}
=T_{d_2/d}\{W\{K\}\} 
$$ that all the points of
the form 
$$\cos (\alpha+\frac{2\pi d_2}{d_1}j), \ \ \ j=1,2, ... ,\frac{d_1}{d}-1,$$
belong to $\sigma_2^{-1}\{K_2\}.$ Since the numbers $d_2/d$ and $d_1/d$ are coprime this implies that 
for any $k=1,2, ... , (d_1/d)-1$ the number
$\cos (\alpha+\frac{2\pi d}{d_1}k)$ belongs to $\sigma_2^{-1}\{K_2\}.$

\vskip 0.3cm
\noindent{\it Remark.} Instead of the condition \be \la{mnb} {\rm card}\{K\}\geq \LCM\{d_1,d_2\} \ee in the formulation of the theorem one can require that \be \la{ner} {\rm card}\{K_1\}\geq d_2/d+1 \ \ \ \ \  {\rm or} \ \ \ \ \  
{\rm card}\{K_2\}\geq d_1/d+1.\ee Indeed,  
for any polynomial
$f(z)$ and any finite set $K\subset \C$ we have:
$$ {\rm card} \{f^{-1}\{K\}\}\geq 
\deg f(z)\,{\rm card} \{K\}-\deg f^{\prime}(z)=$$
\be \la{in} 
=\deg f(z)({\rm card} \{K\}-1)+1.\ee Therefore, any of inequalities \eqref{ner} implies 
inequality \eqref{mnb}.

\vskip 0.3cm
\noindent{\it Proof of corollary 1.} It is enough to observe that 
if ${\rm card}\{K_1\}\geq 2$ then \eqref{in}
implies that 
$${\rm card}\{K_1\}\geq \deg f_1(z)+1=\LCM(d_1,d_2)+1.$$

\vskip 0.3cm
\noindent{\it Proof of corollary 2.} Let $C$ be the circle of the smallest radius containing
the set $K$ and $c$ be its center. Observe that any  
$\mu(z)\in\Sigma_{K}$ transforms $C$ to itself. 
Therefore, $\Sigma_{K}$ is a subgroup
of the group $S^1.$ Since $K$ is a compact set, it follows that if $\Sigma_{K}$
is infinite then $K$ contains with a point $x$ all the circle with the center $c$ containing $x$ and hence is a union of circles. Moreover, setting 
$$\tilde f(z)=z^{d_1} \circ \sigma_1, \ \ \ \ \tilde K_1= \tilde f \{K\}$$ 
where $\sigma_1(z)=z-c$ 
we see that then
$$\tilde f^{-1}\{ \tilde K_1\}=f^{-1}\{K_1\}=K.$$ It follows now from corollary 1 that
$f(z)=\sigma_2 \circ \tilde f(z)$ for some linear function $\sigma_2(z)$ and hence 
$f(z)=\sigma_2 \circ z^{d_1} \circ \sigma_1.$

Suppose now that the group $\Sigma_{K}$ 
is finite. Without loss of generality we can suppose that $c=0.$ Then
$\Sigma_{K}$ is generated by $\varepsilon_b=exp(2\pi i/b).$ Since 
$$f^{-1}\{K_1\}=(f\circ \varepsilon_b z)^{-1}\{K_1\}=K,$$
it follows from corollary 1 that 
\be \la{tyu} f(\varepsilon_bz)=\mu\circ f(z)\ee
for some $\mu(z)\in \Sigma_{K_1}.$ If $\mu(z)=
\alpha z+\beta,$ $\alpha, \beta\in \C,$ then  
$\alpha=\varepsilon_b^a,$
where $a$ is the remainder after division of $\deg f (z)$ by $b$. This implies in particular the equality 
$f(0)=(f(0)-\beta)/\varepsilon_b^a.$
Consider now
the rational function $g(z)=(f(z)-f(0))/z^a.$ 
Since
$$g(\varepsilon_b z)=\frac{f(\varepsilon_b z)-f(0)}{\varepsilon_b^a z^a}=\frac{f(z)-(f(0)-\beta)/\varepsilon_b^a}{z^a}=
g(z)
$$ it is easy to see that $g(z)$ has the form $R(z^b)$
for some polynomial $R(z)$ and hence 
$f(z)=\sigma \circ z^aR(z^b),$ where $\sigma(z)=z+f(0).$

\section{Solutions of $f_1^{-1}\{T\}=f_2^{-1}\{T\}=K$} 
{\it Proof of theorem 2.} If $d_1$ is a divisor of $d_2$ then 
the theorem follows from theorem 1
so we may concentrate on the case when 
$d_1$ is not a divisor of $d_2.$ Let us suppose additionally that $d_1/d>2;$ the case when 
$d_1/d=2$ will be considered separately.

Since $d_1$ is not a divisor of $d_2$ theorem 1 implies that 
there exist non-linear polynomials  
$\tilde f_1(z),$ $\tilde f_2(z),$ $g_1(z),$ $g_2(z),$ 
$\deg g_1(z)=d_2/d,$ $\deg g_2(z)=d_1/d,$ and a polynomial $W(z),$ $\deg W(z)=d,$
satisfying
\eqref{2}, \eqref{kom}. Since in course of the proof of the theorem we will repeatedly use theorem 1 to uniform the notation set 
$$L_1=T, \ \ \ A_1(z)=\tilde f_1(z), \ \ \ B_1(z)=\tilde f_2(z).$$
It follows from \eqref{2}, \eqref{kom} that
$$ A_1^{-1}\{L_1\}= B_1^{-1}\{L_1\}=W\{K\}.$$ 

Furthermore,
by theorem 1 there exists a compact subset of $\C,$ which we denote by $L_2,$
such that  
\be \la{ux} g_1^{-1}\{L_2\}=g_2^{-1}\{L_2\}=L_1\ee 
and there exist linear functions $\sigma_1(z), \sigma_2(z)$ such that
either \be \la{u1} A_1(z)=\sigma_1 \circ z^{d_1/d}, \ \ \ 
B_1(z)=\sigma_2 \circ z^cR(z^{d_1/d}),\ee or 
\be \la{u2}
A_1(z)=\sigma_1 \circ T_{d_1/d}, \ \ \ 
B_1(z)=\sigma_2 \circ T_{d_2/d}\ee
holds.

Set now
\be \la{rys1} 
A_2(z)=g_2(\tilde f_1(z)), \ \ \ \ \ B_2(z)=g_1(\tilde f_2(z)). \ee Then the conditions of the theorem imply
that \be \la{kot} 
A_2^{-1}\{L_2\}=B_2^{-1}\{L_2\}=W\{K\}, \ee
where $\deg A_1(z)=(d_1/d)^2,$ $\deg B_2(z)=(d_2/d)^2.$ Furthermore,
applying theorem 1 to equality \eqref{ux} we conclude that there 
exist polynomials $h_1(z), h_2(z),$ $\deg h_1(z)=d_1/d,$ $\deg h_2(z)=d_2/d.$ such that 
$$h_1(g_1(z))=h_2(g_2(z))$$ and
$$h_1^{-1}\{L_3\}=h_2^{-1}\{L_3\}=L_2$$ for some compact set $L_3\subset \C.$ 
Then for polynomials 
$$A_3(z)=h_1 \circ g_2 \circ f_1, \ \ \  B_3(z)=h_2 \circ g_1 \circ f_2$$
we have: $$ A_3^{-1}\{L_3\}=B_3^{-1}\{L_3\}=W\{K\}, $$
where $\deg A_3(z)=(d_1/d)^3,$ $\deg B_3(z)=(d_2/d)^3.$

Continuing in the same way we conclude that 
for any $r\geq 0$ there exist a compact set $L_r$ and polynomials
$A_r(z), B_r(z),$ $\deg A_r(z)=(d_1/d)^r,$ $\deg B_2(z)=(d_2/d)^r,$
such that \be \la{rav} A_r^{-1}\{L_r\}= B_r^{-1}\{L_r\}=W\{K\},\ee
where for polynomials $A_1(z), B_1(z)$ either \eqref{u1} or \eqref{u2} holds. 
Furthermore, applying theorem 1 to equality \eqref{rav} for $r\geq 2$ we
see that there exist linear functions 
$\sigma_{r,1}(z),$ $\sigma_{r,2}(z),$ $\omega_r (z)$ such that either
\be \la{11} A_r(z)=\sigma_{r,1} \circ z^{(d_1/d)^r} \circ  \omega_r, 
\ \ \ B_r(z)=\sigma_{r,2} \circ
z^{c_r}R_r(z^{(d_1/d)^r}) \circ \omega_r \ee
for some $R_r(z)$ and $c_r,$ 
or
\be \la{22} A_r(z)=\sigma_{r,1} \circ T_{(d_1/d)^r} \circ  \omega_r, 
\ \ \ B_r(z)=\sigma_{r,2} \circ
T_{(d_2/d)^r} \circ \omega_r. \ee

Show that if \eqref{u1} (resp. \eqref{u2}) holds then \eqref{11} (resp. \eqref{22}) holds for 
all $r\geq 2.$
Consider first the case when \eqref{u2} holds and
show that equality \eqref{11} can not be realized. 
Indeed, observe that the formula $T_n(\cos x)=\cos(n\, x)$ implies that 
$T_n^{\prime}(z)=0$ if and only if $z =\cos (\pi k/n),$ where $k=1,2, ... , n-1.$
In particular, since 
$d_1/d>2,$ the polynomial $T_{d_1/d}(z)$ has at least 
two critical points. It follows now from 
the chain rule that 
the polynomial $A_r(z),$ which is by construction a polynomial in $T_{d_1/d}(z),$
also has at least two critical points. On the other hand, 
the polynomial $\sigma_{1,r} \circ z^{(d_1/d)^r} \circ \omega_r$ has only one critical point. 

Similarly, if formula \eqref{u1} holds then \eqref{22} can not be realized since
$A_r(z)$ is a polynomial in $z^{d_1/d}$ and therefore 
has at least one critical point of the multiplicity $>2$ while
the multiplicity of any critical point of the polynomial $\sigma_{1,r} \circ T_{(d_1/d)^r} \circ \omega_r$ is 2.

Consider now the cases when \eqref{11} or \eqref{22} holds separately.
Suppose first that \eqref{22} holds.
Show at the beginning that for any $r\geq 2$ the equality 
\be \la{plm}
\omega_r(z)=\pm z
\ee holds.
Indeed, we have:
$$A_r(z)=(\sigma_{r,1} \circ T_{(d_1/d)^{r-1}}) \circ  
(T_{d_1/d}\circ \omega_r). $$
Therefore, setting 
$$U_r=(\sigma_{r,1} \circ T_{(d_1/d)^{r-1}})^{-1}\{L_r\},$$
we see that 
$$(T_{d_1/d}\circ \omega_r)^{-1}\{U_r\}=(\sigma_1 \circ T_{d_1/d})
^{-1}\{L_1\}.$$ 
By corollary 1, this implies that 
$$T_{d_1/d}\circ \omega_r=\delta\circ \sigma_1 \circ T_{d_1/d}$$
for some linear function $\delta(z).$ Since both parts of this
equality should have the same critical points it follows easily that
\eqref{plm} holds. 

Furthermore, since $T_n(\pm z)=\pm T_n(z)$
equality \eqref{plm} implies that 
$$ A_r(z)=\tilde \sigma_{r,1} \circ T_{(d_1/d)^r}, 
\ \ \ B_r(z)=\tilde \sigma_{r,2} \circ
T_{(d_2/d)^r},$$ for linear functions 
$\tilde \sigma_{r,1}=\pm \sigma_{r,1},$ $\tilde \sigma_{r,2}
=\pm\sigma_{r,2}.$
In particular, setting $M_1=\sigma_{1}^{-1}\{L_1\}$ and
$M_r = \tilde \sigma_{r,1}^{-1}\{L_r\}$ for $r\geq 2,$
we see that for any $r\geq 1$ the equality 
\be \la{rra} (T_{(d_1/d)^{r}})^{-1}\{M_r\}=
W\{K\}\ee
holds.

The equality \eqref{rra} implies that
the compact set $W\{K\}$ together
with a point $u$ contains all the points $v$ such that 
\be \la{gop} T_{(d_1/d)^{r}}(v)=T_{(d_1/d)^{r}}(u) \ee for some $r\geq 1.$ 
Choose $\alpha \in \C$ such that $u=\cos \alpha.$ Then condition \eqref{gop} is
equivalent to the condition that $W\{K\}$ contains all the points 
of the form
$$\cos \left(\alpha+2\pi\left(\frac{d}{d_1}\right)^r j\right), 
\ \ \ j=1,2, ...  ,\left(\frac{d_1}{d}\right)^r-1,$$
where $r\geq 1.$ Since $W\{K\}$ is a compact set  
it follows that $W\{K\}$ contains all the set $E_{\alpha }=\cos(\alpha+s),$ $0\leq s \leq 2\pi.$
It is easy to see that $E_{\alpha }$ is an ellipse
which in the coordinates $x=\Re z, y=\Im z$ is defined by the equation  
$$\frac{x^2}{a^2}+\frac{y^2}{b^2}=1, \ \ \ a=\frac{1}{2}\left(\vert e^{i\alpha}\vert + \frac{1}{\vert e^{i\alpha}\vert}\right),
\ \ \ b=\frac{1}{2}\left(\vert e^{i\alpha}\vert - \frac{1}{\vert e^{i\alpha}\vert}\right).$$

Therefore, we can represent $W\{K\}$ as a union of ellipses
$$W\{K\}=\bigcup_{t\in U} E_{t}$$
for some compact subset $U$ of the segment $[0,i\infty).$
Furthermore, since $T_n\{E_ t \}=
E_{tn}$ we have: \be \la{el1} T=
\sigma_1\{T_{d_1/d}\{W\{K\}\}\}=\bigcup_{t\in U} \sigma_1\{E_{td_1/d}\}.\ee On the other hand, \be \la{el2}T=
\sigma_2\{T_{d_2/d}\{W\{K\}\}\}=\bigcup_{t\in U} \sigma_2\{E_{{td_2/d}}\}.\ee

Denote by $t_1$ the point of $U$
with the maximal modulus.
Then formulas \eqref{el1}, \eqref{el2} imply that 
the ellipses $\sigma_1\{E_{{t_1d_1/d}}\}$
and $\sigma_2\{E_{{t_1d_2/d}}\}$ coincide. In particular they
have the same 
focuses. Since focuses of all ellipses $E_{\alpha}, \alpha\in \C$ are $\pm 1$
we conclude that $\sigma_2\circ \sigma_1^{-1}=\pm z$ and hence \eqref{tt3} holds. 
Furthermore, the equality 
$$\sigma_1\{E_{ {t_1d_1/d}}\}
=\sigma_1\{\pm\, E_{ {t_1d_2/d}}\}$$ implies that
$$\frac{1}{2}\left(\vert e^{it_1}\vert^{d_1/d} + \frac{1}{\vert e^{it_1}\vert^{d_1/d}}\right)=
\frac{1}{2}\left(\vert e^{it_1}\vert^{d_2/d} + \frac{1}{\vert e^{it_1}\vert^{d_2/d}}\right)
.$$
Since $d_2\neq d_1$ this follows $t_1=0.$
Therefore, $W\{K\}=[-1,1]$ and hence 
$$T=\sigma_1\{T_{d_1/d}\{[-1,1]\}\}=\sigma_1\{[-1,1]\}$$
is a segment. 

Consider now the case when \eqref{11} holds.
Since
$$A_r(z)=(\sigma_{r,1} \circ z^{(d_1/d)^{r-1}}) \circ  
(z^{d_1/d}\circ \omega_r), $$
setting 
$$U_r=(\sigma_{r,1} \circ z^{(d_1/d)^{r-1}})^{-1}\{L_r\},$$
we see that 
$$(z^{d_1/d}\circ \omega_r)^{-1}\{U_r\}=(\sigma_1 \circ z^{d_1/d})
^{-1}\{L_1\}.$$ 
By corollary 1, this implies that 
$$z^{d_1/d}\circ \omega_r=\delta\circ \sigma_1 \circ z^{d_1/d}$$
for some linear function $\delta(z).$ Comparing critical points of the both sides of this equality we conclude that
$\omega_r(z)=\gamma_r z$ for some $\gamma_r\in \C.$

Therefore, for $r\geq 2$ we have:
$$A_r(z)=\tilde \sigma_{r,1} \circ z^{(d_1/d)^r}, 
\ \ \ B_r(z)=\tilde \sigma_{r,2} \circ
z^{c_r}\tilde R_r(z^{(d_1/d)^r}) $$
for some linear functions $\tilde \sigma_{r,1},$ $\tilde \sigma_{r,2}$
and a polynomial $\tilde R_r(z).$
In particular, setting $M_1=\sigma_{1}^{-1}\{L_1\}$ and 
$M_r = \tilde \sigma_{r,1}^{-1}\{L_r\}$ for $r\geq 2$
we see that for any $r\geq 1$ the equality 
\be \la{rra1} (z^{(d_1/d)^{r}})^{-1}\{M_r\}=
W\{K\}\ee
holds.

Equality \eqref{rra1} implies
that $W\{K\}$ together with a point $u$ contains 
all the points of the form $\varepsilon u,$ where
$\varepsilon^{(d_1/d)^{r}}=1$ for some $r\geq 0$
and therefore 
all the circle $x^2+y^2 = \vert u \vert.$ It follows that $W\{K\}$ is
a union of such circles and hence by corollary 2 
the function $\tilde f_2(z)$ actually has the form $\sigma_2 \circ z^{d_2/d}$
for some linear function $\sigma_2(z).$ Furthermore,
equality 
\be \la{yt}
T=\tilde f_1(W\{K\})=\tilde f_2(W\{K\}) \ee
implies that 
$\sigma_2(z)=\sigma_1\circ \gamma z$ for some $\gamma\in \C$ and hence
\eqref{tt2} holds.

To finish the proof we only must consider the case when $\deg \tilde f_1(z)=2.$
Define $A_2(z), B_2(z)$ by
formula \eqref{rys1}. 
Since $\deg A_2(z)=4>2$ equality
\eqref{kot} implies that
$L_2$ is either a union of circles or a segment and, respectively, $W\{K\}$ is either a union of circles
centered at the origin or a segment $[-1,1].$ 
If $W\{K\}$ is a union of circles then by corollary 2 we have:
$$\tilde f_1(z)=\sigma_1 \circ z^{d_1/d}, \ \ \ \tilde f_2(z)=\sigma_2 \circ z^{d_2/d}$$
and as above equality \eqref{yt} implies that \eqref{tt2} holds.

On the other hand, if $W\{K\}$ is a segment $[-1,1]$ then the equality 
$$\tilde f_1^{-1}\{T\}=T_2^{-1}\{[-1,1]\}$$ holds and applying corollary 1
we conclude that there exists a linear function $\sigma_1$ such that
$$T=\sigma_1\{[-1,1]\}, \ \ \ \tilde f_1(z)=\sigma_1 \circ T_2.$$
Hence, $T$ is a segment.
Similarly,
$$\tilde f_2^{-1}\{T\}=T_{d_2/d}^{-1}\{[-1,1]\}$$ and 
$$T=\sigma_2\{[-1,1]\}, \ \ \ \tilde f_2(z)=\sigma_2 \circ T_{d_2/d}.$$
It follows now from $\sigma_1\{[-1,1]\}=\sigma_2\{[-1,1]\}$ that
$\sigma_1^{-1}\circ\sigma_2=\pm z$ and hence \eqref{tt3}
holds.

\section{Solutions of $f_1^{-1}\{T\}=f_2^{-1}\{T\}=T$}

\vskip 0.2cm
\noindent{\it Proof of theorem 3.} First of all 
consider the case when
$T$ is finite. In this case inequality \eqref{in} implies that $T$ is a point, $T=t\in \C$. Let $a_1$ be 
the leading coefficient of $f_1(z).$ Set $\sigma(z)=\alpha (z-t),$ where $\alpha^{d_1-1}=a_1.$ Then the polynomial $f(z)=\sigma\circ f_1 \circ \sigma^{-1}$ has the leading coefficient $1$ and satisfies $f^{-1}\{0\}=0.$
It follows that $\sigma\circ f_1 \circ \sigma^{-1}=z^{d_1}.$ Similarly, 
$\sigma\circ f_2 \circ \sigma^{-1}=\gamma z^{d_2}$ for some
$\gamma\in\C.$ 

Assume now that $T$ is infinite.
Consider from the beginning the case when 
the set $T$ is either a union of circles with the common center or a segment. Suppose first that $T$ 
is a union of circles with the common center $c.$ Without loss of generality we can assume
that $c=0.$ Then it follows from corollary 2 that
\be \la{vbn} f_1(z)=\gamma_1 z^{d_1}, 
\ \ \ 
f_2(z)=\gamma_2 z^{d_2}\ee for some $\gamma_1, \gamma_2\in \C.$ Therefore, for 
$\sigma(z) =\alpha z,$ where $\alpha$ satisfies $\alpha^{d_1-1}=\gamma_1,$ the equalities
\eqref{po1} hold with $\gamma=\gamma_2\alpha^{-d_2+1}.$
Furthermore, if $r=\max_{z\in T}\vert z\vert $ then \eqref{vbn} and \eqref{xx}
imply that 
$\gamma_1 r^{d_1-1}=\gamma_2 r^{d_2-1}.$ Therefore, since $r>0,$ the equality 
\be\la{yhn} 
\gamma_1^{d_2-1}=\gamma_2^{d_1-1}\ee
 holds
and hence 
\be \la{nhy} \vert \gamma\vert =\vert \gamma_2 \vert \vert \alpha^{-d_2+1} \vert=
\vert \gamma_2 \vert \vert \gamma_1^{-\frac{d_2-1}{d_1-1}} \vert=1.\ee

Similarly, if
$T$ is a segment then setting
$$p_1(z)= \sigma \circ T_{d_1}\circ \sigma^{-1}, \ \ \ p_2(z)=\sigma \circ T_{d_2} \circ \sigma^{-1}
,$$ where $\sigma(z)$ is a linear function such that $T= \sigma\{[-1,1]\}$ 
we see that
$$p_1^{-1}\{T\}=p_2^{-1}\{T\}=T.$$ By the corollary 1
this implies that
that $$f_1(z)=\delta_1 \circ p_1(z), \ \ \ 
f_2(z)=\delta_2 \circ p_2(z), 
$$ for some linear function $\delta_1(z),$ $\delta_2(z)$
such that $\delta_1\{T\}=\delta_2\{T\}=T.$
Since any linear function which transforms $T$ to itself
has the form 
$$\delta(z)=\sigma \circ \pm z \circ \sigma^{-1}$$
it follows that \eqref{po2} holds.

Consider now the case when $T$ is 
distinct from a union of circles or a segment.
Let $p(z),$ $s=\deg  p(z),$ be a non-linear polynomial of the minimal 
degree satisfying
\be \la{z} f^{-1}\{T\}=\{T\}. \ee
Show that then for any polynomial
$q(z),$ $t=\deg q(z),$ 
satisfying \eqref{z} the equality 
$t=s^k$ holds for some $k\geq 1.$ Indeed, suppose that  
$s^k < t  < s^{k+1}$ for some $k\geq 1.$ Since
$$\left(p^{\circ k}(z)\right)^{-1}\{T\}=q^{-1}\{T\}=T$$
it follows from theorem 2 that there exists a polynomial 
$r(z)$ such that $$q(z)=r( p(z)), \ \ \ r^{-1}\{T\}=T.$$
Since $1 <\deg r(z) < s$ this contradicts to the assumption
about $p(z).$ 

Therefore, $\deg f_1(z)=s^{k_1},$ $\deg f_2(z)=s^{k_2}$
for some $k_1,k_2\geq 1.$ Since 
$$\left( p^{\circ k_1}(z)\right)^{-1}\{T\}=f_1^{-1}\{T\}, \ \ \ 
\left( p^{\circ k_2}(z)\right)^{-1}\{T\}
=f_2^{-1}\{T\}$$ it follows now
from corollary 1 that equalities \eqref{po3} hold. 
Furthermore, since $T$ is distinct 
from a union of circles corollary 2 implies that 
$\Sigma_T$ is finite.

\vskip 0.2cm
\noindent{\it Proof of theorem 4.} Prove at first the equivalence of conditions 1 
and 2. Clearly, it is enough to show that any of the conclusions $a),b),c)$
in the formulation of theorem 3 implies that polynomials $f_1(z),$ $f_2(z)$
have the same Julia sets. In the cases $a),b)$ this is obvious so suppose that the case $c)$ holds. 
Denote by $J_{f_1}$ and $J_{f_2}$ the Julia sets of the polynomials $f_1(z),$ $f_2(z)$
and by $K_{f_1}$ and $K_{f_2}$ their filled-in Julia sets. 
Since for any polynomial $f(z)$ the equality $J_{f}=\partial K_{f}$ holds
in order to prove the equality $J_{f_1}=J_{f_2}$ 
it is enough to prove that $K_{f_1} =K_{f_2}.$

Without loss of generality we can assume that the center of the disk of the smallest radius
containing $T$ is zero. Then 
\be f_1(z)= \eta_1 p^{\circ s_1} , \ \ \ f_2(z)=\theta_2  p^{\circ s_2}
\ee for some $b$-th roots of unity $\eta_1, \theta_1.$ 
Furthermore, applying formula \eqref{tyu} to the polynomial $p(z)$ and taking into account that
$\mu(z)=\alpha z$ for some $\alpha\in \C$ we see that $p(z)=z^aR(z^b)$ for some polynomial $R(z).$ 
This implies that for any $j\geq 1$ we have: 
$$f_1^{\circ j}(z) =\eta_j p^{\circ s_1 j}, \ \ \ \  
f_2^{\circ j}(z) =\theta_j p^{\circ s_2 j}$$
for some $b$-th roots of unity $\eta_j, \theta_j.$ Therefore, the equalities 
$$\vert f_1^{\circ j}(z) \vert =\vert  p^{\circ s_1 j}\vert , \ \ \ \  
\vert f_2^{\circ j}(z) \vert =\vert  p^{\circ s_2 j}\vert $$ hold. This follows that
$K_{f_1}=K_{f_2}=K_{p}$ and hence  
$J_{f_1}=J_{f_2}=J_p.$

Prove now the equivalence of conditions 2 and 3.
Suppose first that 2 holds and set 
$J=J_{f_1}=J_{f_2}.$ Then we have:
$$(f_1\circ f_2)^{-1}\{J\}=(f_2\circ f_1)^{-1}\{J\}.$$ It follows
now from corollary 1 that \eqref{urr} holds with $\mu \in \Sigma_J=J_{f_1}\cap J_{f_2}.$

Furthermore, if
\eqref{urr} holds and $\mu(z)=z$ in other words if $f_1(z), f_2(z)$ commute then
the equality $J=J_{f_1}=J_{f_2}$ was already established by Julia \cite{j}
(for any rational functions) and can be proved easily as follows (\cite{a1}). 
Suppose that $z\in K_{f_1}.$ Since \eqref{urr} implies that
$$f_1(f_2^{\circ k}(z))=f_2^{\circ k}(f_1(z))$$ we conclude that
$f_1\{K_{f_2}\}\subset\{K_{f_2}\}$ for any $k\geq 1.$ Hence, 
$f_1^{\circ  j}\{K_{f_2}\}\subset\{K_{f_2}\}$ for any $j\geq 1$ and therefore 
$K_{f_2}\subset K_{f_1}.$ By the symmetry also $K_{f_1}\subset
K_{f_2}$ and hence $K_{f_1}=K_{f_2},$ $J_{f_1}= J_{f_2}.$ 

Consider now the case when $\mu(z)\neq z.$ 
If
$\mu^{\circ j}(z)\neq z$ for any $j$ then both $\Sigma_{T_1}$ and 
$\Sigma_{T_2}$ are infinite and by corollary 2 taking into account that 
$\Sigma_{T_1}\cap \Sigma_{T_2}\neq \varnothing$ we conclude that $\Sigma_{T_1}$ and 
$\Sigma_{T_2}$ are  
unions of circles with the common center $c.$
Furthermore, since \eqref{urr} implies that for any linear function $\nu(z)$
the equality 
$$\tilde f_1(\tilde f_2(z))=\tilde \mu(\tilde f_2(\tilde f_1(z)))$$ holds 
with 
$$\tilde f_1 =\nu \circ f_1 \circ \nu^{-1}, \ \ \ 
\tilde f_2 =\nu \circ f_2 \circ \nu^{-1}, \ \ \ 
\tilde \mu =\nu \circ \mu \circ \nu^{-1},$$
without loss of generality we can assume that 
$c=0.$ 
In this case the corollary 2 implies there exist $\gamma_1, 
\gamma_2 \in \C$ such that equalities \eqref{vbn} hold. 
Therefore, setting $\sigma(z) =\alpha z,$ where $\alpha^{d_1-1}=\gamma_1,$
we see that equalities
\eqref{po1} hold with $\gamma=\gamma_2\alpha^{-d_2+1}.$
Moreover, equality \eqref{urr} implies equality \eqref{yhn} and therefore equality \eqref{nhy}.
Hence, 
$J_{f_1}=J_{f_2}.$

Suppose now that $\mu(z)$ is a rotation of  finite order $d$ around a
point $c.$ As above we may assume that $c=0.$ Then $\mu(z)=\varepsilon_d z$ for some primitive $d$-th root of unity $\varepsilon_d.$ 
Show that  
\be \la{nn}
f_1(z)= z^{a_1}R_1(z^d), 
\ \ \ \  f_2(z)=z^{a_2}R_2(z^d)\ee
for some polynomials $R_1(z), R_2(z)$ and integers $a_1,a_2.$ Indeed,
if both $\Sigma_{T_1}$ and $\Sigma_{T_2}$ are finite 
then, taking into account the equality $c=0,$ we conclude as above that
there exist polynomials $\tilde R_1(z), \tilde R_2(z)$ and integers $\tilde a_1,\tilde a_2$
such that 
$$f_1(z)=z^{\tilde a_1} \tilde R_1(z^{d_1}), \ \ \  
f_2(z)=z^{\tilde a_2} \tilde R_2(z^{d_2}),$$
where $d_i$ is the order of $\Sigma_{T_i},$ $i=1,2$. Since $d\vert d_1, d\vert d_2$
this implies that \eqref{nn} holds.
On the other hand, if one of (or both) groups 
$\Sigma_{T_1},$ $\Sigma_{T_2}$ 
is infinite then 
it follows from \eqref{vbn} that \eqref{nn} holds.

Following \cite{b} define the polynomials  
$$\tilde f_1(z)= z^{a_1}R_1^d(z), 
\ \ \ \  \tilde f_2(z)=z^{a_2}R_2^d(z)$$ and show that they commute.
Indeed, clearly \be \la{jul}\tilde f_i \circ z^d = z^d \circ f_i, 
\ \ \ i=1,2. \ee
Therefore, we have: 
$$\tilde f_1 \circ \tilde f_2 \circ z^d=\tilde f_1 \circ z^d \circ f_2=
z^d \circ f_1 \circ f_2= z^d \circ \mu \circ f_2\circ f_1=$$
$$=z^d \circ f_2\circ f_1=
\tilde f_2 \circ z^d \circ f_1
=\tilde f_2 \circ \tilde f_1 \circ z^d.$$ Hence,
$\tilde f_1 \circ \tilde f_2=\tilde f_2 \circ \tilde f_1.$ 

Furthermore, since $\tilde f_1(z)$ and $\tilde f_2(z)$ commute 
we have $K_{\tilde f_1}=K_{\tilde f_2}.$ On the other hand, \eqref{jul}
implies that  $$\tilde f_i^{\circ j} \circ z^d = z^d \circ
f_i^{\circ j}, \ \ \ i=1,2 $$ for any $j\geq 1.$ This follows that 
$$K_{f_i}=(z^d)^{-1}\{\tilde K_{f_i}\}, \ \ \ i=1,2.$$ Hence,
$K_{f_1}=K_{f_2}$ and $J_{f_1}=J_{f_2}.$

\bibliographystyle{amsplain}

\begin{thebibliography}{10}


\bibitem {a1} P. Atela, {\it Sharing a Julia set: the polynomial case,}  Progress in holomorphic dynamics,  102--115, Pitman Res. Notes Math. Ser., 387, Longman, Harlow, 1998.
 	
\bibitem{a2} P. Atela, J. Hu, {\it Commuting polynomials and polynomials with same Julia set,} Internat. J. Bifur. Chaos Appl. Sci. Engrg. 6 (1996), no. 12A, 2427--2432


\bibitem {b} I. Baker, A. Eremenko, {\it A problem on Julia sets,}  Ann. Acad. Sci. Fennicae (series A.I. Math.) 12 (1987), 229--236. 


\bibitem {be1} A. Beardon, {\it Polynomials with identical Julia sets,}
Complex Variables, Theory Appl. 17, No.3-4, 195-200 (1992).

\bibitem {be2} A. Beardon, {\it Symmetries of Julia sets,} 
Bull. Lond. Math. Soc. 22, No.6, 576-582 (1990).


\bibitem {d} T. Dinh, {\it Ensembles d'unicit\'e pour les polyn\^omes,} Ergodic Theory Dynam. Systems  22  (2002),  no. 1, 171--186. 

\bibitem {d2} T. Dinh, {\it Distribution des pr\'eimages et des points
p\'eriodiques d'une correspondance polynomiale,} Bull. Soc. Math. France
133 (2005), no. 3, 363--394.

\bibitem {er}  A. Eremenko, {\it
On some functional equations connected with iteration of rational functions,} 
Leningr. Math. J. 1, (1990), No.4, 905-919. 

\bibitem {en} H. Engstrom, {\it
Polynomial substitutions,} Amer. J. Math. 63, (1941), 249-255. 

\bibitem {fa} P. Fatou, {\it Sur l'it\'eration analytique et les substitutions permutables,} J. de Math. 2 (1923), 343.


\bibitem {f} J. Fernandez, {\it
A note on the Julia set of polynomials,} Complex Variables, Theory Appl. 12, No.1-4, 83-85 (1989).


\bibitem {j} G. Julia, {\it M\'emoire sur la permutabilit\'e des fractions rationnelles,} Ann. Ecole Norm. Sup. 39 (1922), 131-215.

\bibitem{peh} B. Fischer, F. Peherstorfer, {\it 
Chebyshev approximation via polynomial mappings and the convergence behavior of Krylov subspace methods}, 
Electron. Trans. Numer. Anal. 12 (2001), 205--215 (electronic).

\bibitem {kamo} S. Kamo, P. Borodin, {\it 
Chebyshev polynomials for Julia sets,} Moscow Univ. Math. Bull. 49, no. 5, 44--45 (1995).

\bibitem{c} G. Lorentz, {\it Approximation of functions,} 
Holt, Rinehart and Winston, New-York-Chicago, III.-Toronto, Ont. 1966.

\bibitem {p2} I. Ostrovskii, F. Pakovitch, M. Zaidenberg, {\it A remark on complex polynomials of least deviation,} Internat. Math. Res. Notices (1996), no. 14, 699--703. 

\bibitem {p1} F. Pakovitch, {\it Sur un probl\`eme d'unicit\'e pour les fonctions m\'eromorphes} C. R. Acad. Sci. Paris Sér. I Math.  323  (1996),  no. 7, 745--748. 
 

\bibitem {ri} J. Ritt, \textit{Prime and composite polynomials,} Trans.
Amer. Math. Soc.  23, no. 1, 51--66 (1922).

\bibitem {r2} J. Ritt, {\it Permutable rational functions,} Trans. Amer. Math. Soc. 25 (1923), 399-448.
 

\bibitem {sch} A. Schinzel, \textit{Polynomials with special regard to
reducibility}, Encyclopedia of Mathematics and Its Applications
77, Cambridge University Press, 2000.

\bibitem {sh} W. Schmidt, N. Steinmetz, {\it The polynomials associated with a Julia set,}  Bull. London Math. Soc.  27  (1995),  no. 3, 239--241. 

\bibitem {y} C. Yang, {\it Open problem}, in Complex analysis,
Proceedings of the S.U.N.Y. Brockport Conf. on Complex Function Theory, June 7--9, 1976, Edited by Sanford S. Miller. Lecture Notes in Pure and Applied Mathematics, Vol. 36, Marcel Dekker, Inc., New York-Basel, 1978.






\end{thebibliography}

\end{document}